\documentclass[a4paper,12pt]{article}
\usepackage[dvips]{epsfig}
\usepackage{amsmath,amssymb,euscript,graphicx,amsfonts,pifont}
\usepackage{enumerate,color}
\usepackage{graphicx}
\usepackage{hyperref}
\usepackage{tikz}

\newtheorem{theorem}{Theorem}[section]
\newtheorem{proposition}[theorem]{Proposition}
\newtheorem{lemma}[theorem]{Lemma}

\newtheorem{remark}[theorem]{Remark}

\newcounter{case}

\renewcommand{\thecase}{\arabic{case}}
\newcounter{subcase}

\numberwithin{subcase}{case}

\newcommand{\Qed}{\rule{2.5mm}{3mm}}

\renewcommand{\P}{\mathcal{P}}

\newcommand{\B}{\mathcal{B}}
\renewcommand{\S}{\mathcal{S}}

\newcommand{\ZZ}{\mathbb{Z}}

\newcommand{\FF}{\mathbb{F}}

\def\syl{\hbox{\rm Syl}}
\def\soc{\hbox{\rm Soc}}
\def\PG{\hbox{\rm PG}}
\def\PGL{\hbox{\rm PGL}}
\def\GL{\hbox{\rm GL}}

\def\PSL{\hbox{\rm PSL}}

\def\char{\hbox{\rm char}}
\def\demo{{\bf Proof}\hskip10pt}
\def\di{\bigm|} \def\lg{\langle} \def\rg{\rangle}
\def\qqed{\hfill $\Box$}

\def\a{\alpha} \def\b{\beta}   
  \def\t{\tau}  
\def\th{\theta}

\def\Aut{\hbox{\rm Aut\,}}
  
\def\f{\noindent}

 \def\og{\overline G} \def\oh{\overline H}

   \def\o1{\overline 1}

\def\o{\overline}   
\def\di{\bigm|} \def\lg{\langle} \def\rg{\rangle}

\usepackage{minitoc}

\begin{document}
	\begin{center}
		{\bf\large  Hamilton cycles in vertex-transitive graphs of order $6p$}
		\footnote{This work was supported in part  by the National Natural Science Foundation of China
			(12471332).}
	\end{center}
	
	\begin{center}
		Shaofei Du and Tianlei Zhou\footnote{Corresponding author: 2230501007@cnu.edu.cn}\\
		\medskip
		{\it {\small
				School of Mathematical Sciences, Capital Normal University,\\
				Bejing 100048, People's Republic of China\\}}
	\end{center}
	
	\renewcommand{\thefootnote}{\empty}

\begin{abstract}
	It was shown by Kutnar and \v Sparl in 2009 that every connected vertex-transitive graph of order	$6p$, where $p$ is a prime, contains a Hamilton path. In this paper, it will be shown that every such graph contains a  Hamilton cycle,   except for the Petersen graph by replacing each vertex by a triangle.
\end{abstract}

\section{Introduction}
\label{sec:intro}
\indent
A simple path (resp. cycle) going through all vertices of a graph is called a Hamilton path (resp. cycle).
In 1969, Lov\'asz \cite{L70} asked that if there exists  a
finite  connected vertex-transitive graph of order more than 3 without a Hamilton path; and   in 1981, Alspach  \cite{A81} asked if there exists an infinite number of connected vertex-transitive graphs that do not have a Hamilton cycle.
Till now,  with the exception of $K_2$, only four connected
vertex-transitive graphs that do not have a Hamilton cycle are known to exist.
These four graphs are the Petersen graph, the Coxeter graph and the two graphs obtained
from them by replacing each vertex by a triangle.
The fact that none of these four graphs
is a Cayley graph  has led to a folklore conjecture that except for $K_2$,
every Cayley graph has a Hamilton cycle, see   \cite{CG96,D83,
	GWM14,GM07,GKMM12,KM09,MR234,
	DM83,WM18,WM15}.

The existence of Hamilton paths, and in some cases
also Hamilton cycles, in  connected vertex-transitive graphs
has been shown for graphs of particular orders, such as, $kp$, where $k\le 6$,
$p^j$, where $j\le 5$, and $2p^2$, where throughout this paper $p$  denotes a prime, see \cite{C98,MR4328721,MR4548744,KM08,KMZ12,KS09,DM87,DM92,DM85,MP82,MP83,MS4,MS2,MS5,Z15}
and a survey paper \cite{KM09}.

Recently,
Kutnar, Maru\v si\v c and the first author proved that every connected vertex-transitive graph of order
$pq$, where $p$ and $q$ are primes, has a Hamilton cycle, except for the Petersen graph (see \cite{MR4328721}). As for the connected vertex-transitive graphs of order $2pq$, where $p$ and $q$ are primes,
Tian, Yu and the first author  showed that  except for the Coxeter graph, every graph of such order  contains
a Hamilton cycle, provided the automorphism group acts primitively on its vertices (see \cite{MR4548744}).

In \cite{KS09}, Kutnar and \v Sparl showed that every connected vertex-transitive graph of order
$6p$, where $p$ is a prime, contains a Hamilton path, and in some particular cases also contains a Hamilton cycle.
The main result of this paper is the following theorem.

\begin{theorem}
	\label{the:main}
	Except for the graph obtained from the Petersen graph by replacing each vertex by a triangle, every connected vertex-transitive graph of order $6p$ contains a Hamilton cycle, where $p$ is a prime.
\end{theorem}

This paper is organized as follows: after this introductory section, some preliminaries will be given in Section 2 and  Theorem~\ref{the:main} will be proved in Section 3.
\section{Terminology, notation and some useful results}
\label{sec:pre}
This section consists of four subsections dealing with basic definitions and notation, generalized orbital graphs,
existence of Hamilton cycles in particular graphs, and group- and   finite field-theoretic results, respectively.

\subsection{Basic definitions and notation}
\label{ssec:definition}
\noindent
Throughout this paper, graphs are finite, simple, connected and undirected,
and groups are finite.
Furthermore, a {\em multigraph} is a generalization
of a graph in which multiple edges and loops are allowed.
Given a graph $X$, we denote the vertex set by $V(X)$ and the edge set by $E(X)$, respectively.
For two adjacent vertices $u,v \in V(X)$, we write
$u \sim v$ and denote the corresponding edge by $\{u,v\}$.
Let $U$ and $W$ be two disjoint subsets of $V(X)$. Then
by $X( U )$ and $X[U,W]$ we denote the subgraph of $X$ induced by $U$ and the bipartite subgraph of $X$ induced by the edges having one end-vertex in $U$
and the other end-vertex in $W$, respectively.

Let $G$ be a group acting faithfully and transitively on a set $V$. A nonempty subset $B$ of $V$ is called a $block$ for $G$ if for each $g\in G$ either $B^g=B$ or $B^g\cap B=\emptyset$, where $V$ and the singletons $\{v\}$ $(v\in V)$ are called the {\it trivial} blocks, and other blocks are called {\it nontrivial} blocks. $|B|$ is called the {\it length} of a block $B$. Put $\B=\{B^g: g\in G\}$. Then the sets in $\B$ form a partition of $V$. We call $\B$ the {\em system of blocks} containing $B$.
If $G$ has no nontrivial blocks on $V$, then $G$ is {\em primitive},
and is {\em imprimitive} otherwise. If every nontrivial normal subgroup of $G$ is transitive on $V$, then $G$ is {\em quasiprimitive}, and is {\em non-quasiprimitive} otherwise.
If $G$ acts non-quasiprimitively on $V$ with an intransitive normal subgroup $N$, then for any   $v\in V$, $v^N=\{v^g: g\in N\}$ is called an $N$-{\em block}. The set of $N$-blocks forms a system of blocks for $G$.

A graph $X$ is said to be {\it vertex-transitive} if its automorphism group  $\Aut(X)$ acts transitively on $V(X)$. A vertex-transitive
graph $X$ for which every transitive subgroup of $\Aut(X)$  is primitive is called a {\it primitive graph} and
an {\it imprimitive graph} otherwise.  If $X$ is imprimitive with an imprimitivity block system induced by  a normal subgroup of some transitive subgroup $G$ of $\Aut(X)$, then the graph $X$ is said to be {\it genuinely imprimitive.}  If
$X$ is imprimitive and every transitive subgroup $G\in \Aut(X)$ acts quasiprimitively on $V(X)$, then $X$ is said to be {\it quasiprimitive.}

Let $m\geq 1$ and $n\geq 2$ be integers. An automorphism $\rho$
of a graph $X$ of order $n$ is called $(m,n)$-{\em semiregular}
(in short, {\em semiregular})
if as a permutation on $V(X)$ it has a cycle decomposition consisting
of $m$ cycles of length $n$.
The question whether all vertex-transitive graphs admit a semiregular
automorphism  is one of the famous open problems in algebraic graph theory
(see, for example, \cite{bcc15,seven,DMMN07,G1,M81}).
Let $\P$ be the set of orbits of $\lg \rho\rg$.
Let $A, B \in \P$. By $d(A)$ and $d(A,B)$ we denote the valency of
$X(A)$ and $X[A,B]$, respectively,
reminding  that the graph $X[A,B]$ is
regular. Let the {\em quotient graph corresponding to
	$\P$} be the graph $X_\P$ whose vertex set
is  $\P$ with $A, B \in \P$ adjacent if there
exist vertices $a \in A$ and $b \in B$, such that $a \sim b$ in $X$.
Let the {\em quotient multigraph corresponding to
	$\rho$} be the multigraph $X_\rho$ whose vertex set is $\P$ and in
which $A,B \in \P$ are joined
by $d(A,B)$ edges.
Note that the quotient graph $X_\P$
is precisely the underlying graph of $X_\rho$.

\subsection{Generalized orbital graphs}
\label{ssec:orbital}
A permutation group $G$ on a set $V$ induces
the action of $G$ on $V\times V$, whose
orbits are called {\em orbitals} of $G$, while  $O_0=\{(x,x) \colon x\in V\}$ is said to be {\it trivial}.
Every orbital $O_i$ corresponds to a {\it paired} orbital $O_i^*=\{ (x,y): (y,x)\in O_i\}$, and
$O_i$ is said to be {\it self-paired} if $O_i=O_i^*$.

For an arbitrary union $\cal{O}$ of orbitals (having empty intersection
with $O_0$), the {\em generalized orbital (di)graph} $X(V,\cal{O})$
of the action of $G$ on $V$
with respect to $\cal{O}$ is a simple (di)graph with the vertex set $V$ and
edge set $\cal{O}$.
It can be viewed as an undirected graph if and only if  $\cal{O}$ coincides with its symmetric closure,
that is, $\cal{O}$ has the property that $(x,y)\in \cal{O}$ implies $(y,x)\in \cal{O}$.
Further,  $X(V,\cal{O})$ is
said to be {\em basic orbital graph } if ${\cal O} =O_i\cup O_i^*$ for some $O_i$.

Let $X$ be a graph with a transitive group $G$ of $\Aut(X)$. By $G_v$ and $[G:H]$ we denote the stabilizer of $v\in V(X)$ and the set of cosets of $G$ relative to $H$, respectively, where $H=G_v$.
Then there exists a 1-1 correspondence between orbitals and suborbits of $G$, that is, the orbits of $H$ on $V(X)$. A suborbit corresponding to a {\em self-paired}
orbital is said to be self-paired. Therefore, the (generalized) orbital (di)graph $X(G/H,\cal{O})$ is denoted by $X(G,H,\S)$ too, provided $\cal{O}$ corresponds to a union  $\S$ of suborbits,  while  $\S$ may be replaced by $D$, a union of some double cosets of $H$. In $X(G,H,D)$, $Hg_1\sim Hg_2$ if and only if $Hg_2\in \{Hdg_1:d\in D\}$, where $g_1,g_2\in G$.

Let $X$ be a genuinely imprimitive graph of order $n$ and let $G$ be an imprimitive subgroup of $\Aut(X)$ with a nontrivial intransitive normal $N$. Then $N$ induces blocks of length $m$, where $m\ge 2$ and $m \mid n$. Then $X\cong X(G,H,D)$ for the union of some double cosets of $H$, where $H=G_v$. Let $\B$ be the set of $N$-blocks. Let $A, B \in \cal{B}$. Let the {\em quotient graph corresponding to
	$\cal{B}$} be the graph $X_\B$ whose vertex set
is  $\cal{B}$, and let $A, B \in \cal{B}$. Then $A\sim B$ in $X_\B$ if there
exist two vertices $a \in A$ and $b \in B$, such that $a \sim b$ in $X$. $X_\B$ is also called an {\em $m$-block graph} induced by $N$. Let $\og=G/N$, $\oh=HN/N$ and $\o{D}=\{\overline{Hdh}: d\in D,h\in H\}$. Then $X_\B\cong X(\og, \oh, \o{D})$ implying that $X_\B$ is a vertex-transitive graph. Moreover, $X_\B$ is connected if $X$ is connected.

\subsection{Existence of Hamilton cycles in particular graphs}
\label{ssec:numbers}
The following  known results about existence of Hamilton cycles in particular graphs will be used later.

\begin{proposition}
	\label{pro:2}
	{\rm \cite[Theorem~1]{BJ78}}
	Every $2$-connected regular graph of order $n$ and valency at least $n/3$ contains a Hamilton cycle.
	
\end{proposition}

\begin{proposition}
	\label{pro:4}
	{\rm \cite[Lemma~5]{MP82}}
	Let $X$ be a graph admitting an $(m, p)$-$semiregular$ automorphism $\rho$, where $p$ is a prime. Let $C$ be a cycle of length $k$ in the quotient graph $X_{\mathcal{P}}$, where $\mathcal{P}$
	is the set of orbits of $\rho$. Then the lift of $C$ either contains a cycle of length $kp$ or it consists of $p$ disjoint $k$-cycles. In the latter case, we have $d(S, S'$) = 1 for every
	edge $SS'$ of $C$.
\end{proposition}

\begin{proposition}
	\label{the:main2}
	{\rm \cite[Theorem~1.4]{MR4328721}}
	With the exception of the Petersen graph, every connected vertex-transitive graph of order $pq$, where $p$ and $q$ are primes, contains a Hamilton cycle.
\end{proposition}

\begin{proposition}
	\label{the:main8}
	{\rm \cite[Theorem~1.1]{KS09}}
	Every connected vertex-transitive graph of order $6p$, where $p$ is a prime, contains a Hamilton path. Moreover, with the exception of the truncation of the Petersen graph, every
	such graph which is not genuinely imprimitive contains a Hamilton cycle.
\end{proposition}

\begin{proposition}
	\label{the:main9}
	{\rm \cite[Theorem~1.1]{MR4548744}}
	Except for the Coxeter graph, every connected vertex-transitive graph of order $2pq$ contains a Hamilton cycle provided the automorphism group acts primitively on its vertices, where $p$ and $q$ are primes.
\end{proposition}

\begin{proposition}
	\label{6p} {\rm \cite[Lemma~3.4, 3.5]{KS09}} Suppose that $X$ is a connected vertex-transitive graph of order $6p$, where $p$ is a prime. If a transitive subgroup $G$ of $\Aut (X)$ contains a minimal normal
	subgroup $N$
	which induces  blocks of length  $6$ or $p$,  then $X$ has a Hamilton cycle.
\end{proposition}

\begin{lemma}
	\label{lemma:2}
	Let $X=X(T,H ,D)$ be an orbital graph of order $t$, 
	where $D=HgH\cup Hg^{-1}H$ for some $g\in T\backslash H $.  
	Let $G =T\times \langle c \rangle $, where $|c|=p$ for a prime
	coprime to $t$. If $X$ has a Hamilton cycle, then both  $Y_1=X(G,H,Dc\cup Dc^{-1})$ and $Y_2 = X(G,H,D\cup HcH \cup Hc\sp{-1}H)$ contain a Hamilton cycle, where $k\neq 0$. In particular, if $T$ is primitive on the set $[T:H]$ and $|[G:H]|=2rs$, where $r$ and $s$ are distinct primes, then there exists some $D'=Hg'H\cup Hg'^{-1}H$, where $g'\in T\backslash H$, such that
	the connected orbital graph of $X(G, H, D_0)$ contains a subgraph with the form either $X(T,H,D'c^k\cup D'c^{-k})$ or $X(G,H,D'\cup Hc^kH \cup Hc\sp{-k}H)$.	
\end{lemma}

\demo  Let  $C:\,  H, Hg _ {1}, Hg _ {2}, \ldots , Hg _ {t-1}, H$ be a Hamilton cycle of $X$.

\vskip  3mm

(1) $Y_1=X(G,H,Dc\cup Dc^{-1})$.
\vskip 4mm

Since $Hg_{i}c\sp{j} \sim Hg_{i+1}c\sp{j+1}$ in $Y_1$ and $(t,p)=1$,  we may get a Hamilton cycle of $Y_1$:
$$\begin{array}{ll}  C_1:\, &H,  Hg _ {1}c, Hg _ {2}c\sp{2}, \ldots, Hg_{t-1}c^{t-1}, Hc^t,Hg_1c^{t+1}, \ldots,Hg_{t-1}c^{2t-1},\\
	&Hc^{2t}, Hg_{1}c^{2t+1},\ldots,Hg_{1}c^{(p-1)t+1},\ldots, Hg _ {t-1}c\sp{pt-1}, H. \end{array}$$

\vskip 3mm
(2) $Y_2 = X(G,H,D\cup HcH \cup Hc\sp{-1}H)$.
\vskip 3mm
Note that  there exists the number $p$ of disjoint cycles of length $t$ in $Y_2$:
$$Hc\sp i, Hg _ {1}c\sp i , Hg _ {2}c\sp i,  \ldots , Hg _ {t-1}c\sp i ,Hc\sp i,$$
\f where  $i\in \ZZ_p$. Observing  $Hg_sc\sp i\sim Hg_sc\sp {i\pm 1}$, we get a Hamilton cycle of $Y_2$  as follows:
$$\begin{array}{lll}  C_2:\, &H Hg_{1} \ldots Hg_{t-1}Hg_{t-1}cHg_{t-2}c\ldots  HcHc\sp2Hg_1c\sp2\ldots Hc\sp{-1}H, \,   &{\rm if}\, p=2;\\
	C_3:\, &H, Hc, \ldots, Hc\sp{-1}, Hg_{1}c\sp {-1}, Hg_{1}c\sp {-2}, \ldots,  Hg_{1},&\\
	& Hg_2, Hg_2c, \ldots,  Hg_2,  c\sp {-1}, \ldots,  Hg_{t-1}H, &{\rm if}\, 2\di t;\\
	C_4:\, &H,  Hg_{1}, \ldots  Hg_{t-1}, Hg_{t-1}c,  Hg_{t-2}c, \ldots,  Hg_1c, Hg_1c\sp2, Hg_2c\sp2,\\
	&\ldots,  Hg_{t-1}c^2, Hg_{t-1}c^3, Hg_{t-2}c^3,\ldots, Hg_{1}c^3, \ldots , Hg_{t-1}c\sp{-1},Hc\sp{-1},\\
	& Hc\sp{-2}, \ldots,  Hc, H,  &{\rm if}\, 2\nmid tp. \end{array}$$

Now, suppose that $T$ is primitive. Then every $H$-double coset of $G$ has the form $Hg'c^{i}H$ for some $g'\in T$ and $i\in \ZZ_p$. Clearly, $\lg H, g'c^{i}\rg=G$ if $g'\notin H$ and $i\neq 0$. Suppose that $D'=Hg'H\cup Hg'^{-1}H$ and $D_0'=Hg'c^{i}H\cup Hc^{-i}g'^{-1}H$. Firstly, suppose that $o(c)=2$. Then $Hg'cH= Hg'^{-1}cH$ and so $D_0'=D'c\cup D'c^{-1}$. Secondly, suppose that $o(c)\neq 2$. Then we may assume $T$ is of degree $2r$. Since $T$ is primitive, we know that every suborbits of $T$ is self-paired and $D_0'=D'c^{i}\cup D'c^{-i}$. Thus, if $X(G, H, D_0)$ doesn't contain a subgraph with the form $X(G,H,D'c^k\cup D'c^{-k}))$, then it contains a connected subgraph with the form $X(G,H,Hg_0H\cup Hc^kH \cup Hc\sp{-k}H)$ for some $g_0\in T\backslash H$.
\qqed

\subsection{Group- and  finite field-theoretic results}
\begin{proposition} \label{the:main4} {\rm \cite[Theorem~1.49]{fsg}}
	Every transitive group $G$ of prime degree $p$ has the socle:  $\ZZ_p$;   	$A_p$;   $\PSL(2,11)$; $M_{11}$;   $M_{23}$;  or $\PSL(d,q)$, where $p= (q\sp d-1)/(q-1)$.
\end{proposition}

\begin{proposition}
	\label{2p} {\rm \cite{LS85}}
	Every primitive group $G$ of degree $2p$, where $p$ is a prime, has the socle: $A_{2p}$; $\PSL(2,q)$ where $p=\frac{q+1}2$; $G=M_{22}$;    $A_5$ where $p\in \{3, 5\}$;  Except for $A_5$ of degree 10, the socle of other groups  is $2$-transitive.
\end{proposition}

\begin{proposition} {\rm \cite{LS85}}
	\label{3p}
	Every primitive group $G$ of degree $3p$, where $p\ge 5$ is a prime, has the socle $S$:  $A_{3p}$; $\PSL(d,q)$ where $3p=\frac{q^d-1}{q-1}$, $A_7$ of degree $15$;
	$\PSL(2,7)$ of degree 21; $A_6$ of degree $15$; $A_7$ of degree 21; and  $\PSL(2,19)$ of degree $57$. 	Except for $A_6$ and $A_7$, of degree $15$, the socle $S$  of other groups  has  a metacyclic transitive subgroup $T$ whose normal  subgroup $P$ is normal.
\end{proposition}

\begin{proposition}
	\label{quasi}
	Let $G$ be a quasiprimitive and imprimitive group of degree $pr$, where  $p\ge 5$ is prime and  $r\in \{2, 3\}$. Then  $G$ is an almost simple group with  the socle $T$, where
	
	(1)\,  $r=2$:\,   $T=M_{11}, \PSL(d,q)$, where $d$ is an odd prime, $(d,q-1)=1$ and $p=\frac{q\sp d-1}{q-1}$ (by checking Proposition~\ref{the:main4} directly);
	
	(2)\, $r=3$:\,   $T=\PSL(2,2^{2^s}), \PSL(3,2), \PSL(3,3), \PSL(d,q),$  where $q$ is an odd prime, $(q,r-1)=1$, $3\di (q-1)$ and $p=\frac{q\sp d-1}{q-1}$ (by {\rm \cite{WX}}).
\end{proposition}

\begin{proposition}
	\label{the:main5}
	{\rm \cite[I. Satz~7.8]{fgt}}
	Suppose that $G$ is a finite group with a normal subgroup $H$ and $P$ is a $Sylow$ $p$-subgroup of $H$, then $G=N_G(P)H$.
\end{proposition}

\begin{proposition}
	\label{the:main6}
	{\rm \cite[I. Satz~4.5]{fgt}}
	For a subgroup $H$ of the group $G$, the factor group $N_G(H)/C_G(H)$ is isomorphic to a subgroup of $\Aut(H)$, the group of automorphisms of $H$.
\end{proposition}

\begin{proposition}
	\label{the:main61} {\rm \cite[I. Satz~17.5]{fgt}}
	Let $K$ be a abelian normal subgroup of the group $G$ such that $K\le B\le G$ where $(|K|,|G:B|)=1$. Then $K$ has a complement in $G$ provided $K$
	has a complement in $B$.
\end{proposition}

\begin{proposition}\label{num} {\rm \cite[page 331]{LN1997}}
	Let $f\in \FF_q[x,y]$ be  absolutely irreducible and of degree $d$, and let $N$ be the number of solutions of $f(x,y)=0$ in $\FF_q^2$. Then
	$$|N-q| \leq(d-1)(d-2)q^{\frac12}+d^2.$$
\end{proposition}

\section{The proof of Theorem~\ref{the:main}}
\f {\bf Outline of the proof:}  To prove Theorem~\ref{the:main}, let $X$ be a connected vertex-transitive graph of order $6p$,
where $p$ is a prime. Then we need to consider the following two cases, separately:
\begin{enumerate}
	\item[\rm(1)] {\it    $\Aut(X)$ contains a  transitive subgroup $G$ which contains  a    normal subgroup $N\ne 1$  inducing nontrivial system of blocks on $V(X)$.}
	\item[\rm(2)] {\it Every transitive subgroup $G$ of $\Aut(X)$ acts quasiprimitively on $V(X)$.}
\end{enumerate}

By Propositions~\ref{the:main8} and \ref{the:main9}, $X$ contains a Hamilton cycle provided $X$ is in Case (2) except for the graph obtained from the Petersen graph by replacing each vertex by a triangle. If $p=2$ or $3$, then $X$ has a Hamilton cycle by Theorem 1.1 in \cite{KM08} and Theorem 2.5 in \cite{DM87}, respectively.
From now on, suppose we are in Case (1) and $p\ge 5$. Let $N\ne 1$ be a {\it maximal intransitive normal subgroup} of $G$ on $V(X)$, that is, every nontrivial  normal subgroup of $G$ containing $N$ is transitive on $V(X)$. Let $\B$ be the set of $N$-blocks. Then $G/N$ acts quasiprimitively on $\B$. If $|B|=p$ ($p$ may be 5), then $X$ has a Hamilton cycle by Proposition~\ref{6p}.
The remaining five  cases $|B|\in\{3p,2p,6, 2, 3\}$ will be dealt with in the following  three subsections. In all cases,    Hamilton cycles of $X$ will be found, except for the graph obtained from the Petersen graph by replacing each vertex by a triangle, see Theorems~\ref{th61}, ~\ref{th71},  ~\ref{th53},    ~\ref{th51} and ~\ref{th52}. Thus, Theorem~\ref{the:main} is proved.
\qqed

\subsection{$|B|=3p$}
\label{sec:5}
\begin{theorem}
	\label{th61}
	If $|B|=3p$,   then $X$ has a  Hamilton cycle.
\end{theorem}
\demo Under the hypothesis, one may  set  $\mathcal{B} =\{B_0,B_1\}$ and $N=G_{B_i}$ (the setwise block stabilizer), where $i\in \ZZ_2$. Then $G=N\lg \t\rg$, for some $\t$ of order $2^e$ which swaps $B_0$ and $B_1$.  Totally
we have the following  four cases:
\vskip 3mm
(1) $N$  is unfaithful on some $N$-block;

(2) $N$  is faithful on both $N$-blocks and is  primitive on some $N$-block;

(3) $N$ is faithful, imprimitive  on both $N$-blocks but   is  not quasiprimitive on some $N$-block; and

(4)  $N$ is faithful, imprimitive and  quasiprimitive  on both $N$-blocks.
\vskip 3mm
\f  Then Theorem~\ref{th61} will be proved, if
a  Hamilton cycle  of $X$ is found  in each case, for case (1), see  Lemma~\ref{lemma:51} and Remark~\ref{remark:51};  for case (2), see Lemma~\ref{lemma:52};
for case (3), see  Lemmas~\ref{lemma:53}, \ref{lemma:73} and Remark~\ref{remark:531};  and  for case (4), see  Lemma~\ref{lemma:54}.\qqed

\begin{lemma}
	\label{lemma:51}
	Suppose that $N$  is unfaithful on some $N$-block and $p\neq 5$. Then $X$ has a Hamilton cycle.
\end{lemma}
\demo Suppose that  $N$ is unfaithful on an $N$-block, say $B_0$.  Let $T_0=N_{(B_0)}$ (the pointwise stabilizer) and $T=T_0\times T_0^\t$. Then $T\lhd G$. If $T$ is transitive on both $B_0$ and $B_1$, then
$X[B_0,B_1]\cong K_{3p, 3p}$ and so it has a Hamilton cycle. So suppose that $T$ induces blocks of length $r$, where $r\in \{ 3, p\}$. Let $\B'$ be the set of $T$-blocks. By Proposition~\ref{the:main2}, the $r$-block graph $X_{\B'}$
of order $\frac{6p}r$ (a product of two primes) has a Hamilton cycle.
Take a  Hamilton cycle $C$ of  $X_{\B'}$. Recall that $T$ is unfaithful on some $T$-block. Then there exist two adjacent $T$-blocks $B_i'$ and $B_j'$ in $C$ such that
$d(B_i', B_j')=r$ in $X$. By Proposition~\ref{pro:4}, the cycle $C$ can be lifted to a Hamilton cycle of $X$. \qqed

\begin{remark}\label{remark:51} Suppose that  $p=5$. Then the arguments in Lemma~\ref{lemma:51} are still true,   except for  $r=3$ and the $3$-block graph $X_{\B'}$ is isomorphic to the Petersen graph.   For this case,
	it  will be shown in Lemma~\ref{lemma:73} that $X$ also has a Hamilton cycle.\end{remark}

\begin{lemma}
	\label{lemma:52} Suppose that $N$  is faithful on both $N$-blocks and is  primitive on some $N$-block. Then $X$ has a Hamilton cycle.
\end{lemma}
\demo Under the hypothesis, we may assume that $N$ is primitive on $B_0$. Take $u_0\in B_0$. Since $N_{u_0}$, the stabilizer of $u_0$ in $N$, is maximal in $N$   if and only if $N_{u_0^\t}$ is maximal in $N$, we know that $N$
is primitive on $B_1$.
Set $K =\soc(N)$. Then $K~\char~N\lhd  G$, $K\lhd G$, and so $K$ is transitive on both $B_0$ and $B_1$. Let $G_0=K\lg \tau\rg$ and $K_0=K\lg \t^2 \rg$. Then $G_0$ is  transitive on $V(X)$  and $G_0/K_0\cong \ZZ_2$. Checking all primitive groups of degree $3p$ listed in
Proposition~\ref{3p}, we get that either $K$
contains a metacyclic subgroup  $T$ which is transitive on $B_0$ and contains a normal subgroup $P$, where $P\in \syl_p(T)$; or
$K\cong A_6, A_7$, both  of them are of  degree $15$.

Firstly, suppose that  $K\cong  A_{3p}$,  which is of degree $3p$. Then $K_{u_0}\cong A_{3p-1}$. Since $K$ has two equivalent transitive representations and has two suborbits with length 1 and $3p-1$ on each $K$-block, we get $d(X) \ge3p-1\ge \frac{6p}3$. By Proposition~\ref{pro:2}, $X$  has a Hamilton cycle.

Secondly, suppose that $K\cong L$, where $L \in \{\PSL(d,q),A_7, \PSL(2,7), \PSL(2,19)\}$($A_7$ of degree 21).  Then $B_0=u_0^T=\{u_0^t: t\in T\}$. Let $P\in \syl_p(T)$ (also, $P\in \syl_p(G_0)$). Then $P^\t=P^g$ for some $g\in K$ and $N_{G_0}(P)\ge \lg T,  g\t^{-1}\rg$. Let $u_1=u_0^{g\tau^{-1}}\in B_1$. Then $B_1=B_0^{\tau^{-1}}=B_0^{g\tau^{-1}}=u_0^{Tg\tau^{-1}}$, that is, $N_{G_0}(P)$ is transitive on $V(X)$, whose normal  Sylow-$p$ subgroup $P$  induces blocks of length $p$. By Proposition~\ref{6p}, $X$ has a Hamilton cycle.

Thirdly, suppose that $K\cong A_6$, which is of degree $15$. Since $K$ is a nonabelian simple group, $\Aut(K)\cong P\Gamma L(2,9)$, $\Aut(K)/K\cong \ZZ_2^2$ and $G_0/K$ is cyclic, we get that $G_0=K\times\lg\t\rg\cong A_6\times \ZZ_2$, $G_0\cong S_6$, $G_0\cong \PGL(2,9)$ or $G_0\cong M_{10}$.
In the first case, $\lg \t \rg$ induces blocks of length $2$. Since the $2$-block graph induced by $\lg \t \rg$ is a connected vertex-transitive graph of order $15$ with a Hamilton cycle and $G_0/\lg \t \rg$ is primitive on the set of $\lg \t\rg$-blocks, by Lemma~\ref{lemma:2}, $X$ has a Hamilton cycle. In the last three cases, $X$ is isomorphic to a connected orbital graph $X_0=X(G_0,H,D)$ of order $30$, where $G_0\cong S_6$, $\PGL(2,9)$ or $M_{10}$, and $H\cong S_4$. Checking all cases with Magma, we get that $X$ has a Hamilton cycle.

Finally, suppose that $K\cong A_7$, which is of degree $15$. Since $K$ is a nonabelian simple group, $\Aut(A_7)\cong S_7$ and $S_7/A_7\cong \ZZ_2$, we get that $G_0\cong A_7\times \ZZ_2$ or $G_0\cong S_7$. In the former case, by Lemma~\ref{lemma:2}, $X$ has a Hamilton cycle. In the latter case, checking all cases with Magma, we get that $X$ has a Hamilton cycle. \qqed

\begin{lemma}
	\label{lemma:53}
	Suppose that $p\neq 5$ and  $N$ is faithful, imprimitive  on both $N$-blocks but   is  not quasiprimitive on some $N$-block.  Then $X$ has a Hamilton cycle.
\end{lemma}
\demo Under the hypothesis, let $K\ne 1$ be a maximal intransitive normal subgroup of $N$ on $B_0$. Then $K$ induces blocks of length $p$ or 3 on $B_0$.
\vskip 3mm

(1) Suppose that $K$ induces blocks of length $p$ on $B_0$. Then $N/K\lessapprox S_3$ and $K\lessapprox S_p\sp3$. Recall that $\t$ swaps $B_0$ and $B_1$. Then $K^{\t}$ induces blocks of length $p$ on $B_1$. Let $P\in \syl_p(K)=\syl_p(K^{\t})=\syl_p(G)$. Since $N$ is faithful on both $N$-blocks, $N_N(P)K=N$ and $N_N(P)K^{\t}=N$, we know that
$N_N(P)$ is transitive on both $B_0$ and $B_1$. Since  $N_G(P)N=G$ and $N_N(P)\leq N_G(P)$, we get  that $N_G(P)$ is  transitive on $V(X)$ and contains  a normal Sylow-$p$ subgroup $P$ which
induces blocks of length $p$.  By Proposition~\ref{6p},  $X$ has a Hamilton cycle.

\vskip 3mm
(2)    Suppose that $K$ induces blocks of length $3$ on $B_0$. Then $K\lessapprox S_3^p$ and $KK^\t/K$ is  a normal subgroup of $N/K$ which is of degree $p$. Since $p \nmid |KK^\t/K|$,   we get
$KK^\t/K$ is trivial, that is, $K=K^\t$   and so $K\lhd G$.
Let ${\B'}$ be the set of $K$-blocks. Then the 3-block graph $X_{\B'}$  is a connected vertex-transitive graph of order $2p$, where $p\ge 7$. By Proposition~\ref{the:main2}, $X_{\B'}$ has a Hamilton cycle. Firstly, suppose that $K$ is unfaithful on some $K$-block of length $3$. Then using the arguments as in Lemma~\ref{lemma:51}, every Hamilton cycle in $X_{\B'}$ can be lifted to a Hamilton cycle of $X$. Secondly, suppose that $K$ is faithful on every $K$-block. Then $K\lessapprox S_3$. Take a subgroup $M/K$ of order $p$ in $N/K$. Then $M=K\times P$ (as $p\ge 7$), for some $P\in \syl_p(G)=\syl_p(N)$, and $M$ is transitive on both $B_0$ and $B_1$.
Moreover, $P^\t=P^g$ for some $g\in N$, which implies that $M^{\t g^{-1}}=M$.
Thus, $G_0:=\lg K\times P, \t g^{-1}\rg \le N_G(P)$ and so $G_0$  is transitive on $V(X)$ with a normal subgroup $P$. By Proposition~\ref{6p} again, $X$ has a Hamilton cycle.
\qqed

\begin{remark}\label{remark:531}
	The arguments in Lemma~\ref{lemma:53} are still true  provided either $p=5$ and $K$ induces blocks of length $5$; or $K$ induces blocks of length $3$ but $X_{\B'}$ is not isomorphic to the Petersen graph.
\end{remark}
\begin{lemma}
	\label{lemma:73}
	Suppose that $p=5$ and  $N$ is faithful, imprimitive  on both $N$-blocks but   is  not quasiprimitive on some $N$-block.
	Then $X$ has a Hamilton cycle.
\end{lemma}
\demo By Remarks~\ref{remark:51} and \ref{remark:531}, we only need to consider the case that  $N$ has a normal subgroup $K$  which induces 10 blocks of length 3.
Let
${\B'}=\{B_{ij}:i\in \ZZ_2, j\in \ZZ_5\}$  be the set of $K$-blocks, where $B_{ij}\in B_i$. We only need to consider the case when  the $3$-block graph $X_{\B'}$  is isomorphic to the Petersen graph.
With the same arguments as in (2) of Lemma~\ref{lemma:53}, we get $K\lhd G$,
no loss, let $K$ be  the kernel of  $G$ on ${\B'}$.
Since $N/K$ induces two orbits of length 5 on ${\B'}$ and $X_{\B'}$  is isomorphic to the Petersen graph,
we have  $G/K\cong \ZZ_5\rtimes \ZZ_4$, which is a minimal vertex-transitive subgroup of the Petersen graph.
Let $P\in \syl_3(K)$ and $Q\in \syl_5(G)$.

Suppose that $P$ is faithful on some $P$-blocks. Then $\ZZ_3\cong P\lhd G$, $Q\in C_G(P)$ and so $Q\lhd G$. Moreover, $Q$ induces blocks of length $5$.
By Proposition~\ref{6p}, $X$ has a Hamilton cycle.

Suppose that $P$ is unfaithful on some $P$-blocks. Then using the same argument as in Lemma~\ref{lemma:53}, we get that $G_0:=N_G(P)$ is transitive on $V(X)$, $G_0/(K\cap G_0)\cong \ZZ_5\rtimes \ZZ_4$ and $\ZZ_3^n\cong P\lhd G_0$, where $2\le n\le 5$. Without loss of generality, we may assume that $Q\in \syl_5(G_0)$. If $n=5$, then $P$ contains a normal subgroup of order 3 in $G_0$, coming back to the previous case.
Since $\GL(2, 3)$ and $\GL(3, 3)$ do not contain any element of order 5, we get $Q$ char $C_G(P)\lhd G_0$. By Proposition~\ref{6p}, $X$ has a Hamilton cycle. Now, we consider the case $P\cong \ZZ_3^4$.
Then from the structure of the Petersen graph, we know that if $B_{ij}\sim B_{ij'}$ in $X_{\B'}$, then $X[B_{ij}, B_{ij'}]\cong K_{3,3}$. Thus, the valency of $X$ is at least 7.
Consider the quotient graph $X_{\P}$ induced by $Q$, where $\P$ is the set of orbits of $Q$ and $\P=\{P_{i,j}: i\in \ZZ_2, j\in \ZZ_3\}$. One may assume that $P_{i,j}\subset B_i$ and easily know that there exists a complete matching in $X_\P$, whose edges have one end-vertex in $B_0$
and the other end-vertex in $B_1$, respectively. Thus, $X_{\P}$ contains a 2-connected and 3-regular subgraph of order 6. By Proposition~\ref{pro:2}, $X_\P$ has a Hamilton cycle. Since $d(P_{i,j_1}, P_{i,j_2})\ge 2$ for any two distinct $Q$-orbits $P_{i,j_1}$ and $P_{i,j_2}$, where $i\in \ZZ_2$ and $j_1,j_2\in \ZZ_3$, by Proposition~\ref{pro:4}, $X$ has a Hamilton cycle.
\qqed

\begin{lemma}
	\label{lemma:54}
	Suppose that $N$ is faithful, imprimitive and  quasiprimitive  on both $N$-blocks. Then $X$ has a Hamilton cycle.
\end{lemma}
\demo  Let $K= \soc(N)$.  Since $K ~\char ~N\lhd  G$, we know that $G_0:=K\lg \t \rg$ acts transitively on $V(X)$, for some $\t$ of order $2^l$ which swaps two $N$-blocks. Now $X$ is isomorphic to a bicoset graph of $K$. Since  $N$ is faithful, imprimitive and  quasiprimitive  on both $N$-blocks,  by Proposition~\ref{quasi}, we know that   $K=\PSL(3,2)$;  $\PSL(3,3)$;   $\PSL(2,2^{2^s})$; or $\PSL(d,q)$,  where $p= (q\sp d-1)/(q-1)$, $d$ is an odd prime and $(d,q-1)=1$.   For the first two groups, Magma may find a Hamilton cycle of $X$.
\vskip 3mm
Now, we consider the case   $K=\PSL(2,2^{2^s})$.
Let $\FF_q^*=\lg \th \rg$ and $f(x)=x^2+\th^mx+1$ be an irreducible polynomial over $\FF_q$, for some $m$.
Set
$$\begin{array}{ll}
	&\ell={\left[\begin{array}{cc} 0&1\\1&0\end{array}\right]},\,
	t=\left[\begin{array}{cc} \th&0\\0&\th^{-1}\end{array}\right],\,
	u=\left[\begin{array}{cc} 1&1\\0&1\end{array}\right],\,
	s(a,b)=\left[ \begin{array}{cc} a&b\\b&a+b\th^m\end{array}\right],\\ &\\
	&S=\lg s(a,b)\di a, b\in\FF_q,  a^2+b^2+ab\th^m=1 \rg.\end{array}$$
\f Then  ${\rm o} (t)=q-1,\,  {\rm o} (u)=2,\, t^{\ell}=t^{-1} \, \, {\rm and}\, \, S\cong
\ZZ_{2^{2^s}+1}.$
Remind that $K=\lg \ell, t, u\rg =\PSL(2, 2^{2^s}) $ has only one conjugacy class of subgroups isomorphic to $S$.

Set $q=2^s$. Let $T=K_{\infty}=\lg u, t\rg \cong \ZZ_2^{2^s}\rtimes\ZZ_{q-1}$, the point-stabilizer of
$K$, relative to $\infty$ in the projective line $\PG(1, q)$.
Let $H=\lg u, t^3\rg  \le T$ and  $[K:H]$ be the set of right cosets of $K$ relative to $H$. Then $K$ has only one conjugacy class of subgroups isomorphic to $H$.
We shall prove the conclusion by the following three steps.

\vskip 3mm
{\it Step 1: Determination of suborbits of $K$.}
\vskip 3mm

Let $\a=H\in [K:H]$.
Since the group $T$ has two orbits on $[K:T]$ (as $K$ is 2-transitive on $\PG(1,q)$) and
$T=H\cup Ht\cup Ht^2$, we get that  the group $H$  has totally three single point suborbits $\{\a^{t^i}\}$
and three suborbits $\a^{t^i\ell H}$ of length $q$, where $i\in \ZZ_3$.
Since $t^i \ell$ is an involution, every suborbit $\a^{t^i\ell H}$ is self-paired.

Since $H\cap S\leq T\cap S=1$, the group $S$ acts semiregularly on $[K:H]$.
Since the group $S$ acts regularly on $[K:T]$, by $T=\cup_{i=0}^{2}Ht^i$ again,
we know that the group $S$ has three orbits $\a^{t^iS}$ on $[K:H]$, where $i\in \ZZ_3$,
and all of them are of length $q+1$.

\vskip 3mm
{\it Step 2: Determination of  some orbital graphs of $K$.}
\vskip 3mm

Firstly, consider the orbital graph $Y(i)=X(K,H,Ht^i\ell H)$ with $i\in\ZZ_3$, and we shall show
$d(\a^{t^{j}S}, \a^{t^{k}S})\ge2$ in $Y(i)$, where $j, k\in\ZZ_3$.

Remind that $N_1(\a)$, the neighborhood of $\a$ in $Y(i)$, is $\{ \a^{t^i\ell h}: h\in H\}=\{ Ht^i\ell h: h\in H\}$,
where $H=\lg u, t^3\rg$.
Then
$$\begin{array}{lcl}
	N_1(\a^{t^{j}})&=&\{Ht^i\ell ht^{j}: h\in H\}=\{Ht^{i-j}\ell h_1: h_1\in H\}\\
	&=&\{H{\left[\begin{array}{cc} 0&\th^{i-j}\\\th^{-(i-j)}&x\end{array}\right]}: x\in
	\FF_q\}.\\
\end{array}$$
Clearly, $d(\a^{t^{j}S}, \a^{t^{k}S})=|N_1(\a^{t^{j}})\cap \a^{t^kS}|$, which is the number of
solutions $s(a,b)$ of the equation
\begin{eqnarray}\label{Eq1}
	H{\left[\begin{array}{cc} 0&\th^{i-j}\\\th^{-(i-j)}&x\end{array}\right]}=Ht^ks(a,b),
\end{eqnarray}
that is,
$${\left[\begin{array}{cc} \th^{3r} &x_1\\ 0 &\th^{-3r} \end{array}\right]}
{\left[\begin{array}{cc} 0&\th^{i-j}\\\th^{-(i-j)}&x\end{array}\right]}=t^ks(a,b),$$
for some $r$ and $x_1$,  that is,
$$\left[\begin{array}{cc} \th^{-(i-j)}x_1 &xx_1+\th^{3r+i-j}\\
	\th^{-(3r+i-j)}&\th^{-3r}x\end{array}\right]= \left[ \begin{array}{cc}
	a\th^{k}&b\th^{k}\\b\th^{-k}&(a+b\th^m)\th^{-k}\end{array}\right].$$
Then Eq(\ref{Eq1}) holds if and only if
$$\left\{\begin{array}{lcll}
	\th^{-(i-j)}x_1 &=&a\th^{k}, \, \hskip 2cm &(i)\\
	xx_1+\th^{3r+i-j}&=&b\th^{k},\, &(ii) \\
	\th^{-(3r+i-j)}&=&b\th^{-k},\, &(iii) \\
	\th^{-3r}x&=&(a+b\th^m)\th^{-k},\, &(iv)\\
	a^2+ab\th^m+b^2&=&1,\,  &(v)
\end{array}\right.  $$
that is,
$$\left\{\begin{array}{lcll}
	x_1&=&a\th^{i-j+k},\,\hskip 2cm & (i')  \\
	xx_1+\th^{3r+i-j}&=&b\th^{k},\, & (ii)\\
	b&=&\th^{-(3r+i-j-k)},\, &(iii')\\
	x&=&(a+b\th^m)\th^{3r-k},\, &(iv') \\
	a^2+ab\th^m+b^2&=&1.\, &(v)
\end{array}\right. $$
Inserting $(i')$, $(iii')$ and $(iv')$ to $(ii)$, we obtain
$$a(a+\th^{-(3r+i-j-k)+m})\th^{3r+i-j}+\th^{3r+i-j}=\th^{-(3r+i-j-2k)}. \hskip 1cm (ii')$$
Note that $(ii')$ and $(iii')$ may imply $(v)$.
Set  $y=\th^{-r}$ and $c=\th^{-(i-j-k)}$. Then the equation $(ii')$ becomes
\begin{eqnarray}\label{Eq2}
	a^2+c\th^m ay^3+c^2y^6+1=0.
\end{eqnarray}

Conversely, given any solution $(a,y)$ of Eq(\ref{Eq2}),
$x_1$, $b$ and $x$ are uniquely determined by $(i')$, $(iii')$ and $(iv')$, respectively; clearly,
$(ii')$ and so $(ii)$ holds; and finally $(v)$ holds following  $(ii')$ and $(iii')$ hold.

In summary,  $d(\a^{t^{j}S}, \a^{t^{k}S})\ge 2$  if and only if  Eq(\ref{Eq2}) has at least two solutions.
Let $n$ be the number of solutions of Eq(\ref{Eq2}). In what follows we shall show $n\ge 2$.

Firstly, suppose that  Eq(\ref{Eq2}) is  reducible in some field $\FF\ge \FF_q$.
Then there exist $g(y)$ and $h(y)\in\FF[y]$ such that  $(a+g(y))(a+h(y))=0$. That is, $g(y)+h(y)=c\th^m y^3$ and $g(y)h(y)=c^2y^6+1$. Therefore,
one has to set  $g(y)=u_3y^3+u_2y^2+u_1y+u_0$ and $h(y)=v_3y^3+v_2y^2+v_1y+v_0$, where $u_3+v_3=c\th^m$ and $u_2=v_2,u_1=v_1,u_0=v_0$.
By $g(y)h(y)=c^2y^6+1$, we get
$$\begin{array}{lcl}
	c^2y^6+1&=&g(y)h(y)=(u_3y^3+u_2y^2+u_1y+u_0)(v_3y^3+u_2y^2+u_1y+u_0)\\
	&=&u_3v_3y^6+(u_3+v_3)y^3(u_2y^2+u_1y+u_0)+(u_2y^2+u_1y+u_0)^2\\
	&=&u_3v_3y^6+c\th^m (u_2y^5+u_1y^4+u_0y^3)+u_2^2y^4+u_1^2y^2+u_0^2\\
	&=&u_3v_3y^6+c\th^m u_2y^5+(c\th^mu_1+u_2^2)y^4+c\th^mu_0y^3+u_1^2y^2+u_0^2,
\end{array}$$
which implies that $c\th^mu_0=0$ and $u_0^2=1$, a contradiction.

Secondly,  suppose that  Eq(\ref{Eq2}) is absolutely irreducible. Then it
follows  from Proposition~\ref{num} that
\begin{eqnarray}\label{Eq3}
	n\ge q-20q^{\frac12}-36.
\end{eqnarray}
Suppose $q\ge 484$. Then $q-20q^{\frac12}-36\ge 2$.
Suppose $q\le 483$.  Since $q=2^{2^s}$, $q+1=p$ and $3\di q-1$, we get  $q=16$ or 256.
Checking by Magma,  Eq(\ref{Eq2}) has at least two  solutions for any $c$, over $\FF_{16}$ or $\FF_{256}$.

\vskip 3mm
{\it Step 3: Graph $X$ has a Hamilton cycle.}
\vskip 3mm

Take two vertices $\a $ and $\b$ contained two biparts such that $K_\a=K_\b=H$.
Then $X$ is a bipartite graph,  where one bipart is $\a^{K}$ and the other one is $\b^{K}$.
Note that the group $H$ has 12 suborbits on $V(X)$, saying $\{\a^{t^i}\},\{\b^{t^i}\},\a^{t^i\ell H}$ and $\b^{t^i\ell H}$, where $i\in \ZZ_3$.
Similarly, the group $S$ has 6 suborbits on $V(X)$, saying $\a^{t^iS}$ and $\b^{t^iS}$, where $i\in \ZZ_3$,
and all of them are of length $q+1$. Let $\S$ be the set of orbits of $S$. Since $\lg H,\t,t\rg=\lg T,\t\rg\lneqq G_0$, we get that $X_0=X(G_0,H, D_0)$ is disconnected, where $D_0$ is the union of all orbits of length $1$.
Hence, the quotient graph $X_{\S}$ is also a bipartite graph with an edge either $\{\a^S,\a^{t^i\ell S}\}$ or $\{\a^S,\b^{t^i\ell S}\}$.

Suppose that $\b^{t^i\ell H}\subset N_1(\a)$, the neighborhood of $\a$,  where $i\in \ZZ_3$. Then $X_{\S}$ contains a subgraph isomorphic to $K_{3,3}$ and  $d(B_1, B_2)\ge 2$ for any two adjacent blocks $B_1$ and $B_2$, see Figure 1.
Hence, the following Hamilton cycle of $X_{\S}$ can be lifted to a Hamilton cycle of $X$:
$$\a^{S}, \b^{S}, \a^{tS}, \b^{tS}, \a^{t^2S}, \b^{t^2S}, \a^{S}.$$

\begin{figure}[!h]
	\centering
	\begin{tikzpicture}[every node/.style={circle, draw, minimum size=10mm}]
		\foreach \x/\name in {1/{${\a}^{S}$}, 2/{${\a}^{tS}$}, 3/{${\a}^{t^2S}$}} {
			\node (L\x) at (-2,-\x*2) {\name};
		}
		
		\foreach \x/\name in {4/{${\b}^{S}$}, 5/{${\b}^{tS}$}, 6/{${\b}^{t^2S}$}} {
			\node (R\x) at (2,-\x*2+6) {\name};
		}
		
		\foreach \x in {1,...,3}
		\foreach \y in {4,...,6}
		\draw (L\x) -- (R\y);
	\end{tikzpicture}
	
	\caption{\small \label{fig2} The quotient graph $X_{\S}$  with respect to the $(6,p)$-semiregular automorphism $\rho=s(0,1)$ when $\b^{t^i\ell H}\subset N_1(\a)$.}
\end{figure}
\begin{figure}[htbp]
	\centering
	\begin{tikzpicture}[every node/.style={circle, draw, minimum size=10mm}]
		\node (L1) at (-2,-2){$\a^{S}$};
		\node (L2) at (-1,-4) {$\a^{tS}$};
		\node (L3) at (-2,-6) {$\a^{t^2S}$};
		
		\node (R4) at (2,-2){$\beta^{t^iS}$};
		\node (R5) at (1,-4) {$\beta^{t^{1+i}S}$};
		\node (R6) at (2,-6) {$\beta^{t^{2+i}S}$};
		
		\draw (L1) --  (L2);
		\draw (L2) --  (L3);
		\draw (L3) --  (L1);
		\draw (R4) --  (R5);
		\draw (R5) --  (R6);
		\draw (R6) --  (R4);
		\draw (L1) --  (R4);
		\draw (L2) --  (R5);
		\draw (L3) --  (R6);
	\end{tikzpicture}
	
	\caption{\small \label{fig3} The quotient graph $X_{\S}$ with respect to the $(6,p)$-semiregular automorphism $\rho=s(0,1)$ when $\a^{t^j\ell H}\cup\b^{t^{i}H}\subset N_1(\a)$.}
\end{figure}
Suppose that $\a^{t^j\ell H}\subset N_1(\a)$ with $j\in \ZZ_3$.
Then the induced subgraph of every bipart  contains a subgraph which is isomorphic to $Y(j)$.
Due to the connectivity of $X$, there must be an edge between $\a^{K}$ and $\b^{K}$, saying   $\a ,  \b^{t^i}$, where $i\in \ZZ_3$.
Then  $X_{\S}$ contains a perfect matching between two biparts, see Figure 2.
Hence,  the following Hamilton cycle of $X_{\S}$ can be lifted to a Hamilton cycle of $X$:
$$\a^{S}, \, \a^{tS},\,  \a^{t^{2}S}, \, \b^{t^{2+i}S},\,  \b^{t^{1+i}S},\,  \b^{t^{i}S},\,  \a^{S}.$$
This finish the case $K=\PSL(2,2^{2^s})$.
\vskip 3mm

Finally, let $K=\PSL(d,q)$, where $\frac{q^d-1}{q-1}=p$ is a prime, $d\ge  3$ and $(d, q-1)=1$.  Let $\B_0$ and $\B_1$ be the set of 3-blocks  of $K$ on two $K$-orbits.  Let $X_{\B}$ be the 3-block graph where $\B=\B_0\cup \B_1$. Let $H=K_v$, and let $M$ be a maximal subgroup of $K$ containing $H$, where $v\in B_0$. Then $K$ has three single suborbits $\{H^{l^i}\}$ and one suborbit $HgH$ of length $3(p-1)$ on $[K:H]$, where $l\in M$ and $g\in K\backslash M$. Suppose that $d(B_0)\ge 3(p-1)$. Then by Proposition~\ref{pro:2}, $X$ contains a Hamilton cycle. Suppose that $d(B_0)\leq 2$. Then $X[\B_0, \B_1]$ contains a connected vertex-transitive  subgraph isomorphic
to one of the following: (i) $K_{p, p}-pK_{1,1}$; (ii) the incident graph of points and hyperplane of $\PG(d-1, q)$; or (iii)  the nonincident graph of points and hyperplane of $\PG(d-1, q)$.
In all cases, $X[\B_0, \B_{1}]$  contains a Hamilton cycle, say $C$.

In the first two cases,  for any  two adjacent $3$-blocks $B_{0i}\in \B_0$ and $B_{1j}\in \B_{1}$, where $i,j\in \ZZ_p$, check that  $X[B_{0i}, B_{1j}]\cong K_{3,3}$ and so $C$ can be lifted to a  Hamilton cycle of $X$.

In the third case, a Sylow-$p$ subgroup (Single subgroup) $P$ of $K$ has three orbits of length $p$ on each $K$-orbit, while by $\B_0'$ and $\B_1'$ we denote the set of $P$-orbits on two $K$-orbits, respectively.
Check that $X[\B_0', \B_1']$ is a graph of order 6 and for any $B'\in V(X[\B_0',\B_1'])$, the valency of the vertex $B'$ is at least $2$. Then one may
find a Hamilton cycle for that. Moreover,  $d(B'_0, B'_1)\ge 2$ for any  two adjacent $P$-orbits $B_0'\in  \B_0'$ and $B_1'\in \B_1'$. By Proposition~\ref{pro:4},  $X$ has a Hamilton cycle.
\qqed

\subsection{$|B|=2p$}
\label{sec:6}
\begin{theorem}
	\label{th71}  Suppose $|B|=2p$.   Then $X$ has a Hamilton cycle.
\end{theorem}
\demo  Under the hypothesis, one may set  $\mathcal{B} =\{B_0,B_1, B_2\}$, $N=G_{B_i}$ and $G/N\cong \ZZ_3$. Then $G=N\lg \t\rg$, for some $\t$ of order $3^e$, which is transitive on the three $N$-blocks.
As in Subsection 3.1, we  shall deal with  four cases in the following Lemmas~\ref{lemma:61}-\ref{lemma:63}, separately. In all cases, $X$ has a Hamilton cycle. \qqed

\begin{lemma}
	\label{lemma:61}
	Suppose that $N$  is unfaithful on some $N$-block. Then $X$ has a Hamilton cycle.
\end{lemma}
\demo Suppose that  $N$ is unfaithful on an $N$-block, say $B_0$.  Let $T_0=N_{(B_0)}$ and $T_i=T_0^{\t^i}$, where $i\in \ZZ_3$.
Suppose that $T_0$ is transitive on $B_1$. Then $X[B_0,B_1]\cong K_{2p,2p}$. In this case, $d(X)\ge 4p$ and $X$ has  a Hamilton cycle.

Suppose that $T_0$ induces blocks of length $r$ on a block, say $B_1$, where $r\in \{2, p\}$.
Firstly, suppose $r=p$. Then $d(B_i,B_j)=d(B_0,B_1)\ge p$, where $i\neq j$ and $i,j\in \ZZ_3$. Then $d(X)\ge2p$. By Proposition~\ref{pro:2}, $X$ has a Hamilton cycle.
Secondly, suppose that $r=2$. Then $T_0$ is a 2-group and so the normal 2-subgroup $T=T_0T_1T_2$ of $G$ induces blocks of length 2.
Using the same arguments as in Lemma~\ref{lemma:51}, every Hamilton cycle of this 2-block graph  can be lifted to a Hamilton cycle of $X$.
\qqed

\begin{lemma}
	\label{lemma:62}
	Suppose that $N$  is faithful on all $N$-blocks and is  primitive on some $N$-block. Then $X$ has a Hamilton cycle.
\end{lemma}
\demo Suppose that $N$ is primitive on an $N$-block, say $B_0$. Then it is primitive on the other two $N$-blocks as in Lemma~\ref{lemma:52}.  Set $K =\soc(N)$. Then $K\lhd G$. Let $G_0=K\lg \t \rg$. Then $G_0$ is transitive on $V(X)$. By Proposition~\ref{2p}, $K\cong T$, where $T\in\{ A_{2p}$, $\PSL(2,q)$, $M_{22}$, $A_5\}$, and $K$ has  three equivalent representations on the three $N$-blocks.
For the first three groups, $K$  is  2-transitive,  which implies either $d(B_0)\ge 2p-1$ and $d(X)\ge 2p+1$; or $d(B_i, B_j)\ge 2p-1$ and so $d(X)\ge 4p-2\ge  \frac {6p}3$, where $i\neq j$.
By Proposition~\ref{pro:2}, $X$ has a Hamilton cycle.

For  $K\cong A_5$, we get  $G_0=K\times\lg\t\rg\cong A_5\times
\ZZ_3$. Let $H_0=(G_0)_v\cong S_3$, where $v\in V(X)$.
Then $X$ contains a subgraph $X_0\cong X(G_0, H_0, D)$, where $D$ is a union of some double cosets of $H_0$. Let $G_1=G_0.\ZZ_2\cong\ZZ_{3}\times S_5$, and let $H_1\cong \ZZ_2\times S_3$ be the subgroup of $G_1$ containing $H$. Let $D_1=H_1DH_1$, a union of some double cosets of $H_1$, and let $X_1\cong X(G_1,H_1,D_1)$. Let $f$ be a bijection from $V(X_0)$ to $V(X_1)$ such that $f(Hg)=H_1g$, where $g\in G_0$. Then we find that $Hg_1\sim Hg_2$ in $X_0$ if and only if $H_1g_1\sim H_1g_2$ in $X_1$, where $g_1,g_2\in G_0$, and so $f$ is also an isomorphism from $X_0$ to $X_1$. Therefore, $\ZZ_3\times S_5$ is isomorphic to a transitive subgroup of $\Aut(X_0)$. Moreover, $\Aut(X_0)$ has a transitive subgroup which is isomorphic to $\ZZ_3\times
(\ZZ_5\rtimes \ZZ_4)= \ZZ_5\rtimes \ZZ_{12}$. By Proposition~\ref{6p}, $X$ has a Hamilton cycle.
\qqed

\begin{lemma}
	\label{lemma:64}
	Suppose that $N$ is faithful, imprimitive  on all $N$-blocks but is not quasiprimitive on some $N$-block.   Then  $X$ has a Hamilton cycle.
\end{lemma}
\demo  Under the hypothesis, let $K\ne 1$ be a maximal intransitive normal subgroup of $N$ on $B_0$. Then $K$ induces blocks of length $p$ or $2$ on $B_0$.
\vskip 3mm

(1)	Suppose that $K$ induces blocks of length $p$ on $B_0$.  Then $N/K\cong \ZZ_2$ and $K\lessapprox S_p\sp2$. Let $P\in \syl_p(K)=\syl_p(N)$. Using the same argument as in Lemma~\ref{lemma:53}, we know that $N_G(P)$ is transitive on $V(X)$, whose normal subgroup $P$ induces blocks of length $p$. By Proposition~\ref{6p}, $X$ has a Hamilton cycle.

\vskip 3mm
(2) Suppose that  $K$ induces blocks of length $2$ on $B_0$.  Then $K=\ZZ_2^l$ for some $l$. Since $K^\t\lhd N$ and $(K^\t K)/K\lhd N/K\le S_p$, we get $K^\t=K.$
As before, we only need to consider the case $l=1$, that is, $K=\ZZ_2$. Let $P\in \syl_p(N)=\syl_p(G)$. Then $KP=K\times P$.
Moreover,  $P^\t=P^g$ for some $g\in N$,  that is, $P^{\t g^{-1}}=P$.
Thus, $S:=\lg K\times P,  \t g^{-1}\rg $ is transitive on $V(X)$, whose  normal Sylow $p$-subgroup $P$ induces blocks of length $p$. By Proposition~\ref{6p} again, $X$ has a Hamilton cycle.
\qqed

\begin{lemma}
	\label{lemma:63}
	Suppose that $N$ is faithful, imprimitive and  quasiprimitive  on all $N$-blocks.  Then  $X$ has a Hamilton cycle.
\end{lemma}
\demo Let $K=\soc(N)$. Then by Proposition~\ref{quasi}, either $K\cong M_{11}$ or $K\cong \PSL(d,q)$, where $d$ is an odd prime and $p=(q\sp d-1)/(q-1)$.

Suppose that $K\cong M_{11}$. Then $\Aut(K)=K$ and so $G$ has a transitive subgroup $K\times \lg \tau\rg\cong M_{11}\times \ZZ_3$. Since $K$ contains $3$ suborbits of length $1$, $1$, $20$ on each $N$-block and $X$ is connected, we get that either $d(B_0)\ge 20$ and $d(X)\ge 22$; or $d(B_0,B_1)\ge 20$ and $d(X)\ge 40$. By Proposition~\ref{pro:2}, $X$ has a Hamilton cycle.

Secondly, suppose that $K\cong \PSL(d,q)$. In this case, three representations of $K$ on three $N$-blocks $B_0$, $B_1$ and $B_2$  must be equivalent. Let $H=K_v$, and let $M$ be the maximal subgroup of $K$ containing $H$, where $v\in B_0$. Then $K$ has two single suborbits $\{H^l\}$ and one suborbit of length $2(p-1)$ on $[K:H]$, where $l\in M\backslash H$. Then either $d(B_0)\ge 2(p-1)$ and $d(X)\ge 2(p-1)+2$; or $d(B_0,B_1)\ge 2(p-1)$ and $d(X)\ge 4(p-1)$, by Proposition~\ref{pro:2}, $X$ contains  a Hamilton cycle.
\qqed

\subsection{$|B|\in \{6, 2, 3\}$}
\label{sec:4}

\begin{theorem}
	\label{th53}  Suppose $|B|=6$.   Then $X$ has a Hamilton cycle.
\end{theorem}
\demo Suppose  $|B|=6$. Then we have the following three cases.

Firstly, suppose that   $N$  does not contain any normal subgroup  $N_1$ of $G$ which induces blocks of length $2$ or $3$. Then $X$ has a Hamilton cycle by Proposition~\ref{6p}.

Secondly, suppose $N_0\le N$ induces blocks of length 3. Let $N_1$ be the kernel of $G$ on the set  $\B'$ of $N_0$-blocks. Then $N_1\leq N$ and $G/N_1\lessapprox \ZZ_2^p\rtimes S_p$.  Suppose that $N_1$ is unfaithful on some $N_0$-block. Then in every  Hamilton cycle  $C_{\B'}$ of $X_{\B'}$ (by Proposition~\ref{the:main2}),   there exist two adjacent $N_0$-blocks $B_1', B_2'$ such that $X[B_1',B_2']\cong K_{3,3}$, so that $C_{\B'}$ can be lifted to a Hamilton cycle of $X$.
So we assume $N_1$ acts faithfully,  that is, $\ZZ_3\cong Q \le N_1\le S_3$. Pick up $N_1\le G_1\le G$ such that $G_1/N_1\cong \ZZ_2^l\rtimes \ZZ_p.$   Then $Q\lhd G_1$ and  $C_{G_1}(Q)$ contains a normal subgroup $L\cong \ZZ_2^l$ such that $Q\times L\lhd G_1$. Then $L\lhd G_1$ and it induces blocks of length 2.  If $l\ne 1$, then $L$ is unfaithful on $L$-blocks so that $X$ has a Hamilton cycle as before. If $l=1$, then the transitive subgroup $(Q\times L).P $ is abelian, where $P\in \syl_p(G_1),$ so that  $P$ induces blocks of length $p$ and then $X$ has a Hamilton cycle by Proposition~\ref{6p}.

Finally, suppose $N_1\le N$ induces blocks of length $2$. Then we have the completely same
arguments as in the third paragraph.\qqed

\begin{theorem}
	\label{th51}
	Suppose that $|B|=2$. Then  $X$ has a Hamilton cycle.
\end{theorem}
\demo  Under the hypothesis, the 2-block graph $X_\B$ is a connected vertex-transitive graph
of order $3p$. By Proposition~\ref{the:main2}, $X_\B$ has a Hamilton cycle.
Suppose that $N$ is unfaithful on a block in $\B$. Then every
Hamilton cycle $C_\B$ in $X_\B$ can be lifted to a Hamilton cycle of $X$ as in Lemma~\ref{lemma:51}. Therefore, we assume that $N$ is faithful on each $N$-block, that is, $N\cong \ZZ_2$.
Since  $G/N$ is a quasiprimitive (including primitive) group of degree $3p$,
by Propositions~\ref{3p} and ~\ref{quasi},   $L/N=\soc(G/N)$ is a nonabelian simple group.
Since $L_{B_0}=N\times L_v$, where $B_0$ is an $N$-block and $v\in B_0$, we have
$N\le N\times L_v\le L$ such that $(|N|, |L:L_{B_0}|)=(2, 3p)=1$. By Proposition~\ref{the:main61}, $L=N\times T$ , where $T\cong L/N$ is a simple group of degree $3p$.
If $T$ is transitive on $V(X)$,  then $T$ is quasiprimitive on $V(X)$, which has been done in \cite{KS09}. So we assume that $T$ is intransitive on $V(X)$ but induces blocks of length $3p$. Following Theorem~\ref{th61},
$X$ has a Hamilton cycle.
\qqed

\begin{theorem}
	\label{th52}
	If $|B|=3$, then $X$ has a Hamilton cycle.
\end{theorem}
\demo   In this case, $X_\B$ is a connected vertex-transitive graph of order $2p$  and $G/N$ acts faithfully and quasiprimitively on $\B = \{B_i:i\in \ZZ_{2p}\}$, listed in  Propositions~\ref{3p} and ~\ref{quasi}.  Set  $L/N=\soc( G/N)$, a nonabelian simple group.

As in Lemma~\ref{lemma:51}, we just need to consider the case when  $N$ is faithful on each $N$-block, that is, $\ZZ_3\cong Q \le N\lessapprox S_3$.  Let $C=C_L(Q)$.
Clearly,  $Q\lhd L$ and so $C\lhd L$.  By Proposition~\ref{the:main6}, we get
$|L:C|\le 2$ and so $|L/Q: C/Q|\le 2$.
Since $C/Q\cong CN/N\cong L/N$, we get $C\cong Q.(L/N)$ and $C$ is a transitive subgroup of $G$ on $V(X)$.
Similar to Theorem~\ref{th51},  one may get  $C=Q\times T$, where $T\cong L/N$ is a simple group of degree $2p$, and just consider the case when $T$ is intransitive on $V(X)$ but induces blocks of length $2p$. Following Theorem~\ref{th71},
$X$ has a Hamilton cycle.\qqed
\bigskip
\begin{center}
	{\bf Acknowledgment}
\end{center}
\vskip 2mm
The authors would like to thank  referees
for their  helpful comments and suggestions.


{\footnotesize 	\begin{thebibliography}{99}
		\itemsep=0pt
		
		
		
		
		\bibitem{A81} B. Alspach,
		The search for long paths and cycles in vertex-transitive graphs and digraphs, in:
		Kevin L. McAvaney (ed.), {\it Combinatorial Mathematics Vlll}, Lecture Notes in Mathematics {\bf 884},
		Springer-Verlag, Berlin, 1981, 14--22.
		
		
		\bibitem{bcc15}
		P.~J.~Cameron (ed.),
		Problems from the Nineteenth British Combinatorial Conference,
		{\em Discrete Math.} {\bf 167/168} (1997), 605--615.
		
		\bibitem{seven}
		P.~J.~Cameron, M.~Giudici, G.~A.~Jones, W.~M.~Kantor,
		M.~H.~Klin, D.~Maru\v{s}i\v{c}, L.~A.~Nowitz,
		Transitive permutation groups without semiregular subgroups,
		{\em J. London Math. Soc.} {\bf 66} (2002), 325--333.
		
		\bibitem{C98}
		Y. Q. Chen,
		On Hamiltonicity of vertex-transitive graphs and digraphs of order $p\sp 4$,
		{\em J. Combin. Theory Ser. B} {\bf 72} (1998), 110--121.
		
		\bibitem{CG96}
		S. Curran and J. A. Gallian,
		Hamiltonian cycles and paths in Cayley graphs and digraphs - a survey,
		{\em Discrete Math.} {\bf 156} (1996), 1--18.
		
		\bibitem{DMMN07}
		E.~Dobson, A.~Malni\v c, D.~Maru\v si\v c and L.~A.~Nowitz,
		Semiregular automorphisms of vertex-transitive graphs of certain valencies,
		{\em J. Combin. Theory Ser. B} {\bf 97} (2007), 371--380.
		
		
		
		\bibitem{MR4328721}
		F. S. Du, K. Kutnar and D.~Maru\v si\v c,
		Resolving the Hamiltonian problem for vertex-transitive graphs of
		order a product of two primes,
		{\it Combinatorica} {\bf 41} (2021), 507--543.
		
		\bibitem{MR4548744}
		F. S. Du, Y. Tian and H. Yu,
		Hamilton cycles in primitive graphs of order $2rs$,
		{\it Ars Math. Contemp} {\bf 23} (2023), \# P3.05
		
		\bibitem{D83}
		E.~Durnberger,
		Connected Cayley graphs of semidirect products
		of cyclic groups of prime order by abelian groups are hamiltonian,
		{\em Discrete Math.} {\bf  46} (1983),  55--68.
		
		\bibitem{GWM14}
		E. Ghaderpour and D. W. Morris,
		Cayley graphs on nilpotent groups
		with cyclic commutator subgroup are hamiltonian,
		{\em Ars Math. Contemp.} {\bf 7} (2014), 55--72.
		
		\bibitem{G1}
		M.~Giudici,
		Quasiprimitive groups with no fixed point free elements of prime order,
		{\em J. London Math. Soc.} {\bf 67} (2003), 73--84.
		
		\bibitem{GM07}
		H. H. Glover and D. Maru\v si\v c,
		Hamiltonicity of cubic Cayley graph,
		{\em J. Eur. Math. Soc.} {\bf 9} (2007), 775--787.
		
		\bibitem{GKMM12}
		H. H. Glover, K. Kutnar, A. Malni\v c and D. Maru\v si\v c,
		Hamilton cycles in $(2,odd,3)$-Cayley graphs,
		{\em Proc. London Math. Soc.} {\bf 104} (2012), 1171--1197.
		
		\bibitem{fsg}
		D. Gorenstein,
		{\it Finite Simple Groups: An Introduction to Their Classifications,} Springer, New York, 1982.
		
		\bibitem{fgt}
		B. Huppert, {\it Endliche Gruppen. I,}
		Springer Berlin, Heidelberg, 1967.
		
		\bibitem{BJ78}
		B. Jackson,
		Hamilton cycles in regular 2-connected graphs.
		{\em J. Combin. Theory Ser. B.} {\bf 29}  (1980),   27--46.
		
		
		
		
		
		\bibitem{KM09}
		K.~Kutnar and D.~Maru\v si\v c,
		Hamilton cycles and paths in vertex-transitive graphs - Current directions,
		{\em Discrete Math.} {\bf 309} (2009), 5491--5500.
		
		\bibitem{KM08}
		K.~Kutnar and D.~Maru\v si\v c,
		Hamiltonicity of vertex-transitive graphs of order $4p$,
		{\em European J. Combin.}  {\bf 29} (2008), 423--438.
		
		
		\bibitem{KMZ12}
		K. Kutnar, D. Maru\v si\v c and C. Zhang,
		Hamilton paths in vertex-transitive graphs of order $10p$,
		{\em European J. Combin.} {\bf 33} (2012), 1043--1077.
		
		\bibitem{KS09}
		K. Kutnar and P. \v Sparl,
		Hamilton paths and cycles in vertex-transitive graphs of order $6p$,
		{\em Discrete Math.} {\bf 309} (2009), 5444--5460.
		
		\bibitem{LN1997}
		R. Lidl and H. Niederreiter, {\it Finite Fields}, Cambridge University Press,
		Cambridge, 1997.
		
		\bibitem{LS85}
		M.~W.~Liebeck and J.~Saxl,
		Primitive permutation groups
		containing an element of large prime order,
		{\em J. London Math. Soc.} (2) {\bf 31} (1985), 237--249.
		
		\bibitem{L70}
		L.~Lov\' {a}sz,
		{ Combinatorial Structures and Their Applications,}
		ed. R.Guy, H.Hanam, N.Sauer and J.Schonheim, Gordon and Breach, New York, 1970.
		
		\bibitem{MR234}
		F. Maghsoudi,
		On Hamiltonicity of Cayley graphs of order $pqrs$,
		{\em Australas. J. Combin} {\bf 84} (2022), 124--166.
		
		\bibitem{DM83}
		D.~Maru\v si\v c,
		Hamiltonian circuits in Cayley graphs,
		{\em Discrete Math.} {\bf 46} (1983), 49--54.
		
		\bibitem{DM87}
		D.~Maru\v si\v c,
		Hamiltonian cycles in vertex symmetric graphs of order $2p\sp 2$,
		{\em Discrete Math.} {\bf 66} (1987), 169--174.
		
		\bibitem{DM92}
		D.~Maru\v si\v c,
		Hamiltonicity of vertex-transitive $pq$-graphs,
		{\em Fourth Czechoslovakian Symposium on Combin.,
			Graphs and Complexity} (1992), 209--212.
		
		\bibitem{M81}
		D.~Maru\v si\v c,
		On vertex symmetric digraphs,
		{\em Discrete Math.} {\bf 36} (1981), 69--81.
		
		
		\bibitem{DM85}
		D.~Maru\v si\v c,
		Vertex transitive graphs and digraphs of order $p\sp k$,
		Cycles in graphs (Burnaby, B.C., 1982) 115--128,
		{\em Ann. Discrete Math.} {\bf 27}, North-Holland, Amsterdam, 1985.
		
		\bibitem{MP82}
		D.~Maru\v si\v c and T.~D.~Parsons,
		Hamiltonian paths in vertex-symmetric graphs of order $5p$,
		{\em Discrete Math.} {\bf 42} (1982), 227--242.
		
		\bibitem{MP83}
		D.~Maru\v si\v c and T.~D.~Parsons,
		Hamiltonian paths in vertex-symmetric graphs of order $4p$,
		{\em Discrete Math.} {\bf 43} (1983), 91--96.
		
		\bibitem{MS4}
		D.~Maru\v si\v c and R.~Scapellato,
		A class of graphs arising from the action of $\PSL(2,q^2)$ on cosets  of $\PGL(2,q)$,
		{\em Discrete Math.} {\bf 134} (1994), 99--110.
		
		\bibitem{MS2}
		D.~Maru\v si\v c and R.~Scapellato,
		Characterizing vertex-transitive $pq$-graphs with an imprimitive subgroup of automorphisms,
		{\em J.~Graph Theory} {\bf 16} (1992), 375--387.
		
		\bibitem{MS5}
		D.~Maru\v si\v c and R.~Scapellato,
		Classifying vertex-transitive graphs whose order is a product of two primes,
		{\em Combinatorica} {\bf 14} (1994), 187--201.
		
		
		
		\bibitem{WM18}
		D. W. Morris,
		Cayley graphs on groups with commutator
		subgroup of order $2p$ are hamiltonian,
		{\em The Art of Discrete and Applied Math.} {\bf 1} (2018),
		$\#$P1.04.
		
		\bibitem{WM15}
		D. W. Morris,
		Odd-order Cayley graphs with commutator subgroup of order $pq$ are hamiltonian,
		{\em Ars Math. Contemp.} {\bf 8} (2015), 1--28.
		
		\bibitem{WX}
		R.J. Wang and M.Y. Xu, A classification of symmetric graphs of
		order $3p$, {\it J. Combin. Theory Ser. B}, {\bf 58} (1993),
		197--216.
		
		\bibitem{Z15}
		J.~Y.~Zhang,
		Vertex-transitive digraphs of order $p^5$ are Hamiltonian,
		{\em Electronic J. Combin.} {\bf 22} (2015), $\#$P1.76.
		
		
	\end{thebibliography}
}
\end{document}